\documentclass[a4paper,11pt]{article}
\usepackage{amssymb}
\usepackage{cite}

\newcommand{\cL}{{\cal L}}

\newcommand{\cC}{\mathcal{C}}
\newcommand{\cA}{\mathbb{A}}
\newcommand{\cB}{\mathcal{B}}

\newcommand{\cH}{\mathcal{H}}
\newcommand{\cK}{\mathcal{K}}
\newcommand{\cO}{\mathcal{O}}
\newcommand{\cQ}{\mathcal{Q}}
\newcommand{\cG}{\mathcal{G}}
\newcommand{\bX}{\mathbb{X}}
\newcommand{\cR}{\mathcal{R}}
\newcommand{\Ds}{(-\Delta)^s }
\newcommand{\cP}{\mathcal{P}}
\newcommand{\bP}{\mathbb{P}}
\newcommand{\R}{\mathbb{R}}

\newcommand{\cZ}{\mathcal{Z}}
\newcommand{\N}{\mathbb{N}}

\newcommand{\cF}{{\cal F}}
\newcommand{\bS}{\mathbb{S}}
\newcommand{\cuad}{{\sqcap\kern-.68em\sqcup}}

\newcommand{\norm}[1]{\|#1\|}

\newtheorem{theorem}{Theorem}[section]
\newtheorem{proposition}{Proposition}[section]

\newtheorem{lemma}{Lemma}[section]
\newtheorem{corollary}{Corollary}[section]
\newtheorem{remark}{Remark}[section]

\begin{document}

\begin{center}{\bf \large    On fractional harmonic functions

 }\bigskip

 \medskip

  {\small

   {   Huyuan Chen\footnote{chenhuyuan@yeah.net}   \quad  and\quad  Ying Wang\footnote{yingwang00@126.com}
   }
  \bigskip

  Department of Mathematics, Jiangxi Normal University,\\   Nanchang, Jiangxi 330022, PR China  
  } \bigskip

\begin{abstract}
Our concern in this paper is to study the qualitative properties for  harmonic functions related 
to the fractional Laplacian. 
Firstly, we classify the polynomials in the whole space and in the half space  for  the fractional Laplacian defined 
in a principle value sense at infinity.   
Secondly, we study  the fractional harmonic functions in half space 
with singularities on the boundary
and the related distributional identities.     
\end{abstract}

  \end{center}
  
   \smallskip

 \noindent {\small {\bf AMS Subject Classifications}:  35R11; 35B40. 
 
  \smallskip


\vspace{2mm}

\setcounter{equation}{0}
\section{Introduction}
Let   $N\in \N$, $s\in(0,1)$ satisfy $N>2s$, 
and $\displaystyle (-\Delta)^s =\lim_{\epsilon\to0^+}  (-\Delta)^s_\epsilon $ be the fractional Laplacian,
where   
  \begin{equation}\label{def 1.1} 
  (-\Delta)^s_\epsilon  u(x)=  c_{N,s}\int_{\cO_\epsilon  }\frac{ u(x)-
u(x+z)}{|z|^{N+2s}}  dz,
\end{equation}
where $\cO_\epsilon$ is  a sequence of   domains in $\R^N$  with certain symmetric properties such that 
$$\bigcup_{\epsilon\in(0,1]} \cO_\epsilon=\R^N$$
and
$$c_{N,s}=2^{2s}\pi^{-\frac N2}s\frac{\Gamma(\frac{N+2s}2)}{\Gamma(1-s)}>0$$  is a normalized constant such that the Fourier transform holds
 $$\cF((-\Delta)^s  u)(\xi)=|\xi|^{2s} \cF(u)(\xi)\quad{\rm for}\ \, u\in L^1(\R^N)$$
 and $\Gamma$ is the Gamma function.  Here the domains with certain symmetric properties could 
 be formed by the balls, the cubes and polyhedrons. In this article, we mainly take 
 $$\cO_\epsilon=B_{\frac1\epsilon}\setminus B_\epsilon,$$ 
where  $B_r(A)$ is the ball centered at $A$ and with the radius $r$, $B_r=B_r(0)$.

  In the last decades, there has been a renewed and increasing interest in the study of Dirichlet problems involving linear and nonlinear integro-differential operators, especially for the simplest model --- the fractional Laplacian. This growing interest is justified both by  seminal advances in the understanding of nonlocal phenomena from a PDE or a probabilistic point of view, see e.g.
\cite{CS0,CS1,RS,musina-nazarov} and the references therein, and by important applications. 
The most important   $s$-harmonic functions are the fractional  Green kernels, which started as early as in 1961 in \cite{BGR} working on a ball, and  the fractional Poisson kernels in bounded domains in \cite{CLO,NV}. For  the fractional Green kernels we refer the readers to  \cite{BS,K,CS}    in bounded domains,   to \cite{FW} in the half space.  
Furthermore,  we refer to  the  $s$-harmonic functions in cones in \cite{TGV} and  blowing up on the boundary of a regular bounded domain in\cite{A}.

    From the integral-differential form of the fractional Laplacian, a natural restriction for the   function $u$ in (\ref{def 1.1}) is 
 $$\|u\|_{L_s^1(\R^N)}:=\int_{\R^N}\frac{|u(x)|}{1+|x|^{N+2s}} dx<+\infty.$$
Under this restriction, the $s$-harmonic functions in the whole domain or half space can't   grow  over $|x|^{2s}$ at the infinity,   see the references \cite{F,CLO}.   However, the principle value sense at infinity in  (\ref{def 1.1})   allows the fractional Laplacian to work on some special  functions not being  in $L_s^1(\R^N)$ such as the  functions in ${\rm span}\{x_1,\cdots,x_N\}$ and  in  ${\rm span}\{x_ix_j:\, i,j=1,\cdots,N,\, i\not=j\}$,  which is well-defined for the fractional laplacian in $\R^N$ by the  oddness.   Our aim in this article is  to study the $s$-harmonic functions  in the whole space and in upper the half space. 
 
  In the whole space,  let's denote the $m$-order polynomials' set with $k\in\N$, 
  $$\bP_k(\R^N)=\Big\{\sum_{|\alpha|=k}b_\alpha x^\alpha:\, \alpha\in\N^N, b_\alpha\in\R \Big\},$$
  where $\N$ is the set of all nonnegative integers  and $\displaystyle x^\alpha=\prod^N_{i=1}x_i^{\alpha_i}$.
 We first  classify the  $s$-harmonic polynomials  in $\R^N$
  i.e. the classical solution of 
  \begin{equation}\label{eq 1.1-p}
 (-\Delta)^s u=0 \quad    {\rm in }\ \, \R^N. 
 \end{equation}
  For $s\in(0,1]$,  we denote   $\cH^s(\R^N)$   the set of all $s$-harmonic functions in $\R^N$
   and for $m\in\N$ 
\begin{equation}\label{eq 1.1-p-m}
\cH_m^s(\R^N)=\big\{u\in\bP_m(\R^N):\, u\ \,  {\rm is}\  s{\rm -harmonic\ in}\ \, \R^N   \big\}.  \end{equation}
  Particularly,  $\cH_m^1(\R^N)$ is the set of $m$ order  harmonic polynomials in $\R^N$, i.e.
 $$
    -\Delta  u=0 \quad   {\rm in }\ \ \, \R^N. 
  $$
  For $P\in \bP_m$,
 we denote $\displaystyle D_{_P}=\sum_{|\alpha|=m} a_\alpha D^\alpha$ and the inner norm 
 $[P,Q]_m=D_{_P}Q$, then 
 $$\cH^1_m=\cQ_m^{\perp},$$
  where 
 $\cQ_m=\{|x|^2P:P\in \cP_{m-2}\}$ for $m\geq 2$, $\cQ_m=\{0\}$ for $m=0,1$.
 The dimension of $\cH^1_m$ is   $C^{m+N-1}_m-C^{m+N-3}_{m-2}$ for $m\geq 2$.  We refer to the book  \cite{AG} for more properties of $\cH_m(\R^N)$.

 \begin{theorem} \label{th 1.1-p} 
Let $N\geq 2$, $s\in(0,1]$ and $\cH_m^s(\R^N)$ be defined in (\ref{eq 1.1-p-m}).  Then \\
$(i)$ for any $m\in\N$, 
$$\cH_m^s(\R^N)= \cH_m^1(\R^N),$$
$(ii)$  
$$\cH^1(\R^N) \subset \cH^s(\R^N).$$
   \end{theorem}

   Our proof of Theorem \ref{th 1.1-p} is based on the mean value property of 
harmonic functions, which allows us to get  a same harmonic polynomials' set for 
  nonlocal operators with  general kernels in Section 2 under suitable assumptions. From the  mean value property, we obtain that all harmonic functions in $\R^N$ are $s$-harmonic. {\it However it is still  open  if $ \cH^1(\R^N)= \cH^s(\R^N).$} \smallskip

The second part in this paper is devoted to  the $s$-harmonic functions   with isolated singularity on the boundary
\begin{equation}\label{eq 1.1}
 (-\Delta)^s u=0 \quad    {\rm in }\ \, \R^N_+, 
\qquad\ \ u=0\ \ \,  {\rm in }\ \,  \R^N_{*}\setminus\{0\}. 
\end{equation}
 
Note that   the prototype problem
 $$ 
  -\Delta  u=0 \quad    {\rm in }\ \ \, \R^N_+,  \qquad
   u=0  \quad  {\rm on }\ \,  \partial \R^N_+\setminus\{0\}
$$ 
   has   the well-known solution -- the Poisson kernel  
   $$\cP_1(x)=c_N |x|^{-N} x_1\quad {\rm in}\ \,  \R^N_+, $$
   where $c_N>0$ is a normalized constant. 
Moreover, the Poisson kernel  $\cP_1$ is a weak solution of 
 \begin{equation}\label{eq 1.2w}
  -\Delta  u=0 \quad    {\rm in }\ \ \R^N_+,  \qquad  
   u=  \delta_0   \quad  {\rm on }\ \,  \partial \R^N_+  
\end{equation}
in the distributional sense that 
 \begin{equation}\label{id 1.2}
  \int_{\R^N_+} \cP_1(x) (-\Delta) \varphi(x) dx=\frac{\partial  }{\partial x_1 } \varphi(0),\quad \forall \varphi \in C^{1,1}(\R^N_+)\cap C_0^1(\overline{\R}^N_+).
  \end{equation}
 Note that  the  Dirac mass in (\ref{eq 1.2w}) holds in the trace sense that 
 $$\lim_{t\to0^+}\int_{\R^{N-1}} \cP_1(t,x') \zeta(x') dx'= \zeta(0)\quad{\rm for}\ \, \zeta\in C_c(\R^{N-1}).$$
Based on the Green kernel and the Poisson kernel, the semilinear elliptic problem  involving the  Radon measure data $\mu,\, \nu$  
$$ 
-\Delta u+ g(u)=\nu \ \ {\rm in }\ \, \Omega,\qquad  
u=\mu  \ \ {\rm on }\ \, \partial\Omega
$$
has been studied extensively in \cite{Kl,VY} and a survey in \cite{V}. Involving the nonlocal operator,   the semilinear problems with measures could see \cite{CV,Kl,KMG,AP,NV}.

Back to our problem (\ref{eq 1.1}),  we want to study a fundamental solution $\cP_s$ of the fractional Laplacian in the upper half space with the formula:
\begin{equation}\label{e 1.1}
\cP_s(x)=\left\{
 \begin{array}{lll}  
 \cK_s\, |x|^{-N}x_1^s \quad &  {\rm for }\ \, x\in \R^N_+,  \\[2.5mm] 
0  & {\rm for }\ \, x\in \R^N_*\setminus\{0\}, 
\end{array}
 \right.
\end{equation}
 where, in the sequel, we use the following constants 
 $$\cK_s= k_{_N,s} 2^{2s-1} s^{-1}\quad{\rm and}\quad \kappa_{_N,s}=\pi^{-(N/2+1)}\Gamma(N/2)\sin(\pi s). $$
In order to state the distributional identity,  we need the following test functions' space: {\it Given a domain $\Omega$ of $\R^N$, let  $\mathbb{X}_s(\Omega)\subset C(\R^N)$ be the space of functions
$\xi$ satisfying:\smallskip
\begin{enumerate}
\item[(i)] $ \xi$ has the compact support  in  $\bar \Omega$;
\item[(ii)]  $\rho(x)^{-s}\xi$ is continuous  in $\bar \Omega$
and $ \|(-\Delta)^s\xi\|_{L^\infty(\Omega)}<\infty$, where $\rho(x)={\rm dist}(x,\R^N\setminus \Omega)$;
\item[(iii)] there exist $\varphi\in L^1(\Omega, \rho^s dx)$
and $\epsilon_0>0$ such that $|(-\Delta)_\epsilon^s\xi|\le
\varphi$ a.e. in $\R^N_+$, for all
$\epsilon\in(0,\epsilon_0]$.
\end{enumerate} }
 \noindent We remark that $\rho(x)=(x_1)_+$ if $\Omega=\R^N_+$ and $C^\infty_c(\Omega)$ is dense in $\bX_s(\Omega)$, where $t_+= \max\{t,0\}$.

 \begin{theorem} \label{th 1.1} 
Let $\cP_s$ be defined in (\ref{e 1.1}) and $e_1=(1,0,\cdots,0)\in\R^N$. Then  $\cP_s$ is   a solution of (\ref{eq 1.1})  and   verifies that   
\begin{equation}\label{id 1.1} \int_{\R^N_+} \cP_s(x) (-\Delta)^s\varphi(x) dx=\frac{\partial^s }{\partial x_1^s} \varphi(0)\quad{\rm for}\ \,\forall\, \varphi\in \bX_s(\R^N_+),
\end{equation}
where 
 $$
 \frac{\partial^s }{\partial x_1^s} \varphi(0)=\lim_{t\to0^+}\frac{\varphi(te_1)}{t^s}. 
$$
  \end{theorem}

It is remarkable that,  unlike (\ref{eq 1.2w}),    in the trace  sense,   $\cP_s$ doesn't
 have the Dirac mass in the boundary trace, i.e. it isn't a weak solution of 
\begin{equation}\label{eq 1.2-d}
  (-\Delta)^s  u=0 \quad    {\rm in }\ \, \R^N_+, 
\qquad\ \   u=\delta_0   \quad  {\rm on }\ \,     \R^N_{*}, 
\end{equation}
which, from \cite{CHW}, has no weak solutions. Thanks to the nonlocal property of the fractional Laplacian, the source outside of the domain could be taken in the calculation of the first equation of (\ref{eq 1.2-d}) directly, which means  outside source in   (\ref{id 1.1}) 
plays a role  different from (\ref{eq 1.2w}).     

Our next attempt is to understand the  role of outside source in   (\ref{id 1.1}).   To this end, we  approximate the fundamental function $\cP_s$   by the Green kernel and the Poisson kernel in the half space. 
Subject to zero Dirichlet boundary condition in $\R^N_*$, the Green kernel  $\cG_{s,\infty}$
has the formula  (see \cite{FW}): 
\begin{equation}\label{Green 1}
\cG_{s,\infty}(x,y) = \left\{
 \begin{array}{lll}  \frac{k_{_N,s}}{2}|x-y|^{2s-N}\int_0^{\frac{4x_{1}y_{1}}{|x-y|^2}}
 t^{s-1} (t+1)^{-\frac N2}\,dt\quad &{\rm if}\ \, x,y\in \R^N_+,\\[2.5mm] 
0\quad& {\rm if}\ \, x\ {\rm or}\ y\not\in \R^N_+.
\end{array}
 \right.
\end{equation}
It is known that the Poisson kernel in half space  is the following: 
\begin{equation}\label{Poisson 1}
\cP_{s,\infty}(x,y)=\cK_{s}\left(\frac{x_1}{-y_1} \right)^s|x-y|^{-N},\quad\forall\,  x\in \R^N_+,\ y\in \R^N_-,
\end{equation}
where $\R^N_-=(-\infty,0)\times \R^{N-1}.$

\begin{theorem} \label{th 1.2} 
Let $\cP_{s,\infty}$, $\cG_{s,\infty}$ be the  Green kernel and the Poisson kernel in $\R^N_+$  
defined in (\ref{Green 1}) and (\ref{Poisson 1}) respectively, then
$$\epsilon^{-s}\cG_{s,\infty}(\cdot, \epsilon e_1)\to \cP_s \qquad{\rm and}\qquad
 \epsilon^{s}\cP_{s,\infty}(\cdot, -\epsilon e_1)\to \cP_s\quad{\rm as} \ \, \epsilon\to0^+ $$
 in $L_s^1(\R^N_+)$ and uniformly in any compact set of $\R^N_+$.
 \end{theorem}
 
Finally, we   study the solution of 
\begin{equation}\label{eq 1.1-wh} 
 (-\Delta)^s u=0 \quad    {\rm in }\ \  \R^N_+,  
\qquad\ \ u=0  \quad  {\rm in }\ \  \R^N_{-} 
\end{equation}
 subject to the blowing up boundary condition on   $\partial\R^N_+:=\{0\}\times \R^{N-1}$.
 
For a given Radon measure $\mu$ in $\R^N_-$,  denote 
 $$ \cP_{s,\infty}[\mu](x):=\left\{
 \begin{array}{lll}
 \int_{\R^N_*}\cP_{s,\infty}(x,y)  d\mu(y)\quad &  {\rm for }\ \, x\in\R^N_+, \\[2.5mm]
 \mu \quad &  {\rm in }\ \  \R^N_-.
 \end{array}
 \right.
 $$
 Particularly,  let $\mu_\epsilon=\delta_{-\epsilon}(y_1) d\omega_{\R^{N-1}}(y')$ for $\epsilon>0$ and  $d\omega_{\R^{N-1}}$
 is the Hausdorff measure in $\R^{N-1}$ and $\delta_{-\epsilon}(y_1)$ is the Dirac mass in variable $y_1$. For simplicity, we  write it $d\omega_{\R^{N-1}}(y')=dy'$, then we have that for $x\in \R^N_+$
  $$ \cP_{s,\infty}[\mu_\epsilon](x) 
 =\int_{\R^{N-1}}\cP_{s,\infty}\big(x,(-\epsilon,y')  \big)   dy' 
 =\cK_{s}\int_{\R^{N-1}} \left(\frac{x_1}{\epsilon} \right)^s|x-(-\epsilon,y')|^{-N} dy' $$

\begin{theorem} \label{th 1.3} 
Let 
$$\cQ_s(x)=\left\{
 \begin{array}{lll} 
 x_1^{s-1}\quad & {\rm in }\ \  \R^N_{+},\\[2mm] 
 0& {\rm in }\ \  \R^N_{-}.
 \end{array}
 \right.
 $$
 Then  $\cQ_s$ is a solution of (\ref{eq 1.1-wh})
 and it verifies that  
\begin{equation}\label{id 1.3} 
\int_{\R^N_+} \cQ_s(x) (-\Delta)^s\varphi(x) dx=\cC_s\int_{\R^{N-1}} \frac{\partial^s }{\partial x_1^s} \varphi(0,x') dx'\quad{\rm for}\ \,\forall\, \varphi\in \bX_s(\R^N_+),
\end{equation}
where 
$$\cC_s=2 \sqrt{\pi} \, \frac{s}{\sin(\pi s)}\frac{\Gamma(s+\frac12)}{\Gamma(s)}.$$

Moreover, we have that 
$$\cC_s \epsilon^{s}\int_{\R^{N}_*}\cP_{s,\infty}[ \mu_\epsilon]dx\to \cQ_s \quad{\rm}\ \, \epsilon\to0^+ $$
in $L_s^1(\R^N_+)$ and uniformly in any compact set of $\R^N_+$.
   \end{theorem}
   
Note that the $\cQ_s$ is known yet in the classical sense, and  we provide a simple proof, which  is  the Poisson kernel in the half space.  Also, we  can also get the well-known $s$-harmonic function
 $(x_1)_+^s$ and we also refer to  \cite[Theorem 5.1]{CLO} for the uniqueness.

Finally, we try to classify the $s$-harmonic functions  in the half space in (\ref{eq 1.1-wh}) without restriction of the continuity on the boundary.  Denote  
$$\cH^s(\R^N_+):=\big\{u\in C(\R^N_+)  \ {\rm verifies (\ref{eq 1.1-wh}) } \big\}.$$

 \begin{theorem} \label{th 1.2-p} 
 Let $N\geq 3$, $s\in(0,1)$ and
 $$\cA_{s,m}(\R^N_+)=\big\{(x_1)_+^s v(x') :\ v \in \cH^s(\R^{N-1})\big\}\cup \big\{(x_1)_+^{s-1} v (x') :\ v \in \cH^s(\R^{N-1})\big\} $$
 for  $m\in\N$. 

  Then
$$\cA_{s,m}(\R^N_+) \subset \cH^s(\R^N_+)\quad {\rm for}\ \   m\geq 0.$$
  
  \end{theorem}
 
We remark that  $\big\{(x_1)_+^s v(x') :\ v \in \cH^s(\R^{N-1})\big\}\subset C(\R^N)$ and
 since  $(x_1)_+^s$ is an $s$-harmonic function  (\ref{eq 1.1-wh}) continuous up to the boundary.

 \smallskip

  The remainder of this paper is organized as follows.  
     Section 2 is devoted to the $s$-harmonic functions in $\R^N$. 
    In Section 3, we prove
  Theorem \ref{th 1.1} by the Green kernel  $\cG_{s,\infty}^+$ via passing to the limit 
  with $y=\epsilon e_1$ as $\epsilon\to 0$.  Section 4 is devoted to the derivation of the Poisson kernel  $\cP_{s,\infty}$ and then show the approximation to the fundamental solution $\cP_s$ by
  this  Poisson kernel.
 Finally, we obtain the whole boundary explosive solution $x_1^{s-1}$ by  the Poisson kernel $\cP_{s,\infty}$ and give the proof  of Theorem \ref{th 1.3}  and Theorem \ref{th 1.2-p}   in Section 5.
      
 \setcounter{equation}{0}
\section{Harmonic polynomials }

Here we involve general kernel for nonlocal problem
\begin{equation}\label{eq 2.1-p}
\cL_K u=0 \quad \   {\rm in }\ \ \R^N,
 \end{equation}
 where
\begin{equation}\label{def 2.1-p} 
\cL_K  u(x)=\lim_{\epsilon\to0^+}  \int_{B_\frac1\epsilon  \setminus B_\epsilon  }\big( u(x)-
u(x+z)\big) K(|z|) dz
\end{equation}
and the kernel $K:(0,+\infty)\to [0,+\infty)$ is continuous. 
We set 
$$\cH_m^{^K}(\R^N)=\Big\{u\in\bP_m(\R^N): u\ \,  {\rm verifies \  (\ref{eq 2.1-p}) }  \Big\}. $$

 We have the following classifications:
      
 \begin{proposition} \label{th 2.1-p} 
Assume that $N\geq 2$,   $m\in\N$,   $K$ verifies that 
 $$\int_0^{+\infty} K(t) t^{N-1+i} dt\in(0,+\infty]\ \ {\rm for}\  i= 2,3,\cdots,m.$$
Let 
$$\kappa_{i,j}(\epsilon)=\frac{  \int_\epsilon^{\frac1\epsilon} K(t) t^{N-1+j} dt }{
\int_\epsilon^{\frac1\epsilon} K(t) t^{N-1+i} dt},\quad i,j= 2,3,\cdots,m. $$
If   either
\begin{equation}\label{con 2.2-p} 
    \lim_{\epsilon\to0^+}\kappa_{i,j}(\epsilon)=0 \ \ {\rm for\ all}\  i,j= 2,3,\cdots,m,\ i>j
\end{equation}   
or 
\begin{equation}\label{con 2.2-po} 
     \lim_{\epsilon\to0^+}\kappa_{i,j}(\epsilon)=0 \ \ {\rm for\ all}\  i,j= 2,3,\cdots,m,\ i<j,
\end{equation}   
then 
$$\cH_m^{^K}(\R^N)= \cH_m^1(\R^N).$$
 Particularly, for $m=0,1$
$$\cH_m^{^K}(\R^N)= \cH_m^1(\R^N)=\bP_m(\R^N).$$
  \end{proposition}
 \noindent{\bf Proof. }
{\it Part I:   we show $\cH_m^1(\R^N)\subset  \cH_m^{^K}(\R^N)$.
Particularly, for $m=0,1$
$$\cH_m^{^K}(\R^N)= \cH_m^1(\R^N)=\bP_m(\R^N).$$}

   Let $h_m\in \cH_m^1(\R^N)$, then for any $z\in \R^N$, $h_m(\cdot+z)$ is also harmonic in $\R^N$.  
  
Note that  $$h_0(x)=c_0\quad{\rm and} 
 \quad h_1(x+z)=\sum^N_{i=1} c_i (x_i+z_i)$$
 for some $c_i\in\R$ with $i=0,\cdots, N$,
  then 
  $$h_1(x+z) -h_1(x)=\sum^N_{i=1} c_i  z_i.$$ 
  Then for any $\epsilon\in(0,1)$,
  $$\int_{B_\frac1\epsilon  \setminus B_\epsilon  }\big( h_1(x)-
h_1(x+z)\big) K(|z|)   dz=\sum^N_{i=1}  c_i\int_{B_\frac1\epsilon  \setminus B_\epsilon  } z_i K(|z|) dz=0$$
  by the oddness of $h_1$.
  Therefore, we obtain that  for $m=0,1$
  $$\cH_m^{^K}(\R^N)= \bP_m(\R^N) =\cH_m^1(\R^N). $$
  
 {\it Now we deal with the case: $m\geq 2$.}  For $h_m\in \cH_m^1$,  we have that 
   $$(-\Delta)_z\big(h_m(x+z) -h_m(x) \big)=0\quad {\rm for\ any }\ \, z\in \R^N$$
   and the mean value property implies that for any $r>0$
   $$(N|\bS^N |)^{-1}r^{1-N}\int_{\partial B_r}\big(h_m(x+z) -h_m(x) \big)d\omega_r(z)=h_m(x+0) -h_m(x)  =0, $$
 that is,
 $$\int_{\partial B_r}\big(h_m(x+z) -h_m(x) \big)d\omega_r(z)=0.$$
 Therefore,   for any $\epsilon\in(0,1)$,
  \begin{eqnarray*}
 \int_{B_\frac1\epsilon  \setminus B_\epsilon  }\big( h_m(x)-
h_m(x+z)\big) K(|z|)   dz    =   
\int_{\epsilon}^{\frac1\epsilon} K(r)  \int_{\partial B_r}\big( h_m(x)-
h_m(x+z)\big)d\omega_r(z)
  = 0.
   \end{eqnarray*} 
 Now passing to the limit as $\epsilon\to0^+$, we obtain that 
 $h_m\in   \cH_m^{^K}(\R^N)$.  
     
 \medskip

{\it Part II:  we show that 
$$\cH_m^{^K}(\R^N)\subset \cH_m^1(\R^N).$$}

Let 
   $$h_m(x)=\sum_{|\alpha| =m} b_{\alpha} x^\alpha \in \cH_m^{^K}(\R^N),$$ 
   where  $m\geq 2$ and $b_{\alpha}\not=0$.
Then direct computation shows that 
  $$h_m(x+z)=\sum_{|\alpha_1|+|\alpha_2|=m}   u_{\alpha_1}(x) v_{\alpha_2}(z),  $$
     where $\alpha_1,\alpha_2\in \N^m$,    
     $$u_{\alpha_1}\in \bP_{|\alpha_1|}(\R^N)\quad {\rm and}\quad v_{\alpha_2}\in \bP_{|\alpha_2|}(\R^N).$$
 Moreover, 
  \begin{equation}\label{e 2.3-p} 
  h_m(x+z) -h_m(x)=\sum^N_{i=1}    u_{i}(x) z_i+ \sum_{|\alpha_1|+|\alpha_2|=m,  \,  |\alpha_2|\geq 2 }   u_{\alpha_1}(x) v_{\alpha_2}(z), 
    \end{equation}    
where $u_i\in \bP_{m-1}(\R^N)$.
    
  Therefore,  from the definition of fractional laplacian 
 we obtain that  for $m\geq 2$
  \begin{eqnarray}
0 &=&\lim_{\epsilon\to0^+}\int_{B_\frac1\epsilon  \setminus B_\epsilon  }\big( h_m(x)-
h_m(x+z)\big) K(|z|)   dz \nonumber  \\[2mm]& =& \lim_{\epsilon\to0^+} \Big(\sum_{ |\alpha_1|+|\alpha_2|=m,  \,  |\alpha_2|\geq 2  }   u_{\alpha_1}(x) \int_{B_\frac1\epsilon  \setminus B_\epsilon  }   v_{\alpha_2}(z) K(|z|) dz\Big)\nonumber
   \\[2mm]&=&  \lim_{\epsilon\to0^+}\left( \sum_{ |\alpha_1|+|\alpha_2|=m,  \,  |\alpha_2|\geq 2  }   u_{\alpha_1}(x) \int_{\bS^{N}} v_{\alpha_2}(z)  d\omega(z)  \int^{\frac1\epsilon}_ {\epsilon }K(r) r^{|\alpha_2| +N-1} dr   \right)\nonumber
    \\[2mm]&=&  \lim_{\epsilon\to0^+} \sum^m_{j=2 } \left( \int^{\frac1\epsilon}_ {\epsilon }K(r) r^{j +N-1} dr  \Big( \sum_{ |\alpha_1|=m-j, |\alpha_2|=j }   u_{\alpha_1}(x) \int_{\bS^{N}} v_{\alpha_2}(z)  d\omega(z)  \Big)   \right).\label{fund11}
 \end{eqnarray} 
Let 
$$\sigma_j(\epsilon)= \int^{\frac1\epsilon}_\epsilon K(r) r^{ N-1+j} dr, \quad j\in\N $$
and
 $$\cZ_j= \sum_{|\alpha_1|=m-j,|\alpha_2|=j} u_{\alpha_1}(x) \int_{\bS^{N}} v_{\alpha_2}(z)  d\omega(z), $$
then   (\ref{fund11}) could be written as
    \begin{eqnarray}\label{aa-aa}
\lim_{\epsilon\to0^+} \sum_{j=2}^m \cZ_j \sigma_j(\epsilon)   =  0.
 \end{eqnarray}  
Under the assumption  (\ref{con 2.2-p}), we prove that $\cZ_j=0$ for $j=2,\cdots,m$, i.e.
  \begin{equation}\label{e 2.3-p--1}  
  \sum_{|\alpha_1|=m-j,|\alpha_2|=j} u_{\alpha_1}(x) \int_{\bS^{N}} v_{\alpha_2}(z)  d\omega(z) =0.
  \end{equation}
In fact, if $\cZ_m\not=0$, then (\ref{aa-aa}) implies that 
$$\lim_{\epsilon\to0^+} \sigma_m(\epsilon)\Big(\cZ_m+\sum^{m-1}_{j=2}\cZ_j \kappa_{m,j}(\epsilon)\Big)=0  $$
which implies $\cZ_m=0$, thanks to 
$$\lim_{\epsilon\to0^+} \kappa_{m,j}(\epsilon) =0.$$
Inductively, we can obtain $\cZ_j=0$ for $j=2,\cdots,m-1$. 
  
 Under the assumption  (\ref{con 2.2-po}), we show that $\cZ_j=0$ for $j=2,\cdots,m$. 
Indeed, if $\cZ_2\not=0$, then (\ref{aa-aa}) implies that 
$$\lim_{\epsilon\to0^+} \sigma_2(\epsilon)\Big(\cZ_2+\sum^{m}_{j=3} \cZ_j \kappa_{2,j}(\epsilon)\Big)=0,  $$
which implies $\cZ_1=0$, thanks to 
$$\lim_{\epsilon\to0^+} \kappa_{2,j}(\epsilon) =0.$$
Inductively, we can obtain $\cZ_j=0$ for $j=3,\cdots,m$.

 Therefore, we obtain that  for any $r>0$ 
 \begin{eqnarray*} 
&& \int_{\partial B_r}\big(h_m(x+z) -h_m(x) \big)d\omega_r(z)  
\\[2mm]&=&  
 \sum_{ |\alpha_1|+|\alpha_2|=m,  \,  |\alpha_2|\geq 2}   u_{\alpha_1}(x)  \int_{\partial B_r } v_{\alpha_2}(z)  d\omega_r(z) 
 \\[2mm]&=&  
 \sum_{j=2}^m r^{j+N-1}  \Big(\sum_{ |\alpha_1|=m-j, |\alpha_2|=j}     u_{\alpha_1}(x) \int_{\bS^{N}}v_{\alpha_2}(z)  d\omega(z) \Big)
  \\[2mm]&=&0 
  \\[2mm]&=&h_m(x+0) -h_m(x),
 \end{eqnarray*} 
 and from the converse of the mean value property we deduce  that 
 \begin{eqnarray*} 
 -\Delta_x h_m(x)=-\Delta_z  h_m(x+0) =-\Delta_z \big(h_m(x+0) -h_m(x) \big)=0.
 \end{eqnarray*} 
   By the arbitrary of $x$, we obtain 
   $h_m\in \cH_m^1(\R^N)$.  We complete the proof.\hfill$\Box$\medskip
   
   \begin{remark}\label{rem 33-1}
 $(i)$ Let 
  $$\frac1c(1+t)^{-N- \zeta}\leq  K_{1,\zeta}(t)\leq ct^{-N- \zeta}  . $$
  for some $c\geq1$ and $\zeta\leq 2$, 
 then  $K_{1,\zeta}$ satisfies  (\ref{con 2.2-p}). 
  
$(ii)$ Let  $$\frac1c t^{-N- \zeta} \chi_{(0,1)}(t)  \leq   K_{2,\zeta}(t)\leq ct^{-N- \zeta}e^{-t}$$
for $c\geq 1$, where $\chi_{(0,1)}(t)=1$ if $t\in(0,1)$ and  $\chi_{(0,1)}(t)=0$ if $t\geq1$.  

Then for $\zeta\geq m$, $K_{2,\zeta}$ satisfies  (\ref{con 2.2-po}).

   \end{remark}

   \medskip

 \noindent  {\bf Proof of Theorem \ref{th 1.1-p}. }  \noindent  {\it Part $(i)$. } 
  When $K(|z|)=c_{N,s}|z|^{-N-2s}$ with $s\in(0,1)$, then  (\ref{con 2.2-p}) holds ture and $\cL_K=(-\Delta)^s$,  $\cH_m^{^K}(\R^N)=\cH_m^s(\R^N)$, then Theorem \ref{th 1.1-p}
 follows by Proposition \ref{th 2.1-p}.   \smallskip

  {\it Part $(ii)$. } For any  $h \in \cH^1(\R^N)$,  we have that 
   $$(-\Delta)_z\big(h (x+z) -h (x) \big)=0\quad {\rm for\ any }\ \, z\in \R^N$$
   and the mean value property implies that for any $r>0$
   $$(N|\bS^N |)^{-1}r^{1-N}\int_{\partial B_r}\big(h (x+z) -h (x) \big)d\omega_r(z)=h (x+0) -h (x)  =0, $$
then  for any $\epsilon\in(0,1)$,
  \begin{eqnarray*}
 \int_{B_\frac1\epsilon  \setminus B_\epsilon  }\big( h (x)-
h(x+z)\big) K(|z|)   dz    =   
\int_{\epsilon}^{\frac1\epsilon} K(r)  \int_{\partial B_r}\big( h (x)-
h (x+z)\big)d\omega_r(z)
  = 0
   \end{eqnarray*} 
   with $K(t)=t^{-N-2s}$.
 Now passing to the limit as $\epsilon\to0^+$, we obtain that 
 $h \in   \cH ^{^K}(\R^N)$. 
  \hfill$\Box$\medskip

\setcounter{equation}{0}
\section{Approximation by the Green kernel}

Note that  the Green function associated with $\Ds$ in the unit
ball $B_1$  was computed
by Blumenthal, Getoor and  Ray  in \cite{BGR} with the formula 
\begin{eqnarray*}
\cG_{s,1}(x,y) &=& \kappa_{_N,s}  |x-y|^{2s-N} \int_1^{(\psi(x,y)+1)^{1/2}} \frac{(z^2-1)^{s-1}}{z^{N-1}}\,dz \\
&=& \frac{\kappa_{_N,s}}{2}|x-y|^{2s-N}\int_0^{\psi(x,y)}
\frac{z^{s-1}}{(z+1)^{N/2}}\,dz
\end{eqnarray*}
for $x,y\in B_1$ and $\cG_{s,1}(x,y)=0$ if $x \not \in B_1$ or $y \not \in B_1$,
where 
$$\psi(x,y)=\frac{(1-|x|^2)(1-|y|^2)}{|x-y|^2} $$
and $k_{N,s}$ is the normalization constant 
$$\kappa_{_N,s}=\pi^{-(N/2+1)}\Gamma(N/2)\sin(\pi s).$$
  If $N=1=2s$, then direct
computations give
$$
\int_0^{\psi(x,y)}
\frac{z^{-1/2}}{(z+1)^{1/2}}\,dz=2\ln\frac{ 1-xy+ (1-x^2)^{1/2}(1-y^2)^{1/2}}{|x-y|}
$$
and  in this case we have that 
$$
\cG_{s,1}(x,y)=\frac{1}{\pi} \ln\frac{ 1-xy+ (1-x^2)^{1/2}(1-y^2)^{1/2}}{|x-y|}.
$$
By the dilation, the Green function for the ball $B_R$ with $R>0$
is given by
\begin{eqnarray*}
\cG_{s,_R}(x,y) &=& R^{2s-N}\cG_{s,1}\left(\frac{x}{R},\frac{y}{R}\right)
\\[2mm]
&=&\frac{k_{_N,s}}{2}|x-y|^{2s-N}\int_0^{\psi_{_R}(x,y)}
\frac{z^{s-1}}{(z+1)^{N/2}}\,dz
\end{eqnarray*}
with $\psi_{_R}(x,y)=\frac{(R^2-|x|^2)(R^2-|y|^2)}{R^2|x-y|^2}$.  From \cite[(3.1)]{FW}, the  Green function $G_{s,\infty}^+$ in the half space $\R^N_+$ has the formula
\begin{equation}\label{eq:half-space-green-function}
\cG_{s,\infty} (x,y) = \frac{k_{_N,s}}{2}|x-y|^{2s-N}\int_0^{\psi_\infty(x,y)}
\frac{r^{s-1}}{(r+1)^{\frac N2}}\,dr\quad{\rm for}\ \, x,y\in \R^N_+,
\end{equation}
where 
$$\psi_\infty(x,y)=\frac{4x_{1}y_{1}}{|x-y|^2}.$$

Next we obtain $\cP_s$ by passing to the limit of Green kernel $G_{s,\infty}^+(x,y)$.  

\begin{lemma}\label{lm 2.1} 
Let $N\geq 2s$, $\cG_{s,\infty}$ be defined in (\ref{eq:half-space-green-function}), then there holds that
\begin{equation}\label{e 2.0-+}
\epsilon^{-s}\cG_{s,\infty} (\cdot, \epsilon e_1)\to \cP_s\quad  \  {\rm as}\ \, \epsilon\to0^+
\end{equation}
in $L^1_s(\R^N_+)$  and uniformly in any compact set of $\R^N_+$.

Moroever,  for $\varphi\in \bX_s(\R^N_+)$
\begin{equation}\label{e 2.1-00}
 \lim_{\epsilon\to0^+}\int_{\R^N_+}\epsilon^{-s}\cG_{s,\infty} (x, \epsilon e_1)(-\Delta)^s\varphi(x) dx=   \int_{\R^N_+} \cP_s(x) (-\Delta)^s\varphi(x) dx.
  \end{equation}

\end{lemma}
\noindent{\bf Proof. }  We first claim that  for any bounded open set $O$,  $\overline O\subset \overline{\R^N_+}\setminus\{0\}$,  there exists  $\epsilon_0\in(0,\frac18)$ such that 
$$2\epsilon_0e_1\not\in \overline O.$$
Now we want to prove that
\begin{equation}\label{e 2.1-u}
\epsilon^{-s}\cG_{s,\infty} (\cdot,\epsilon e_1)\to  \cP_s \quad{\rm as}\ \ \epsilon\to0^+\quad {\rm uniformly\ in}\ O.  
\end{equation}
 
Note that 
$$|x-\epsilon e_1|^2 = |x|^2\left(1-2 |x|^{-2}x_1\epsilon+|x|^{-2}\epsilon^2\right) $$
and for $x\in \overline O\subset \overline{\R^N_+}\setminus\{0\}$,  $|x|^{-2}$, $|x|^{-2}x_1$ are positive and bounded, taking $\epsilon>0$ small enough,
\begin{eqnarray*} 
\frac{4 x_1  }{|x|^2 } \,\epsilon <  \psi_\infty(x,\epsilon e_1)
 &=&\frac{4 x_1  \,\epsilon }{|x|^2 }\frac{1 }{1-2\epsilon |x|^{-2}x_1+|x|^{-2}\epsilon^2}  
 \\[2mm]&<&\frac{4 x_1  \,\epsilon }{|x|^2 }\frac{1 }{1-2\epsilon |x|^{-2}x_1}
 \\[2mm]&<&\frac{4 x_1  }{|x|^2 }\Big(1+2\epsilon |x|^{-2}x_1 \Big)\,\epsilon,
\end{eqnarray*}
then
\begin{eqnarray*} 
|x|^{2s-N}<|x-\epsilon e_1|^{2s-N}<|x|^{2s-N} \left(1+4(N-2s) |x|^{-2}x_1\epsilon \right),
\end{eqnarray*}
\begin{eqnarray*}
\int_0^{\psi_\infty(x,\epsilon e_1)}
\frac{r^{s-1}}{(r+1)^{\frac N2}}\,dr   <   \int_0^{\psi_\infty(x,\epsilon e_1)} r^{s-1} \,dr  
&=&  s^{-1} \psi_\infty^s(x,\epsilon e_1)
 \\[2mm]
&< & 4^s s^{-1}   \frac{ x_1^s  }{|x|^{2s} } \Big(1+2\epsilon |x|^{-2}x_1 \Big)^s \,\epsilon^s
 \\[2mm]
&<&4^s s^{-1} \frac{  x_1^s  }{|x|^{2s} }\Big(1+2^s\epsilon^s |x|^{-s}  \Big)\epsilon^s
\end{eqnarray*}
and
\begin{eqnarray*}
\int_0^{\psi_\infty(x,\epsilon e_1)}
\frac{r^{s-1}}{(r+1)^{\frac N2}}\,dr  &> & \int_0^{\psi_\infty(x,\epsilon e_1)}\Big( r^{s-1}-Nr^s\Big)\,dr  
\\[2mm]
&=& s^{-1} \psi_\infty(x,\epsilon e_1)^s -\frac{N}{1+s} \psi_\infty(x,\epsilon e_1)^{1+s}
\\[2mm]
&>&4^s s^{-1} \frac{  x_1^s  }{|x|^{2s} } \epsilon^s-\frac{N}{1+s} \frac{4 x_1  }{|x|^2 }\Big(1+2\epsilon |x|^{-2}x_1 \Big)^{1+s}\,\epsilon^{1+s}
\\[2mm]
&>&4^s s^{-1} \frac{  x_1^s  }{|x|^{2s} } \epsilon^s\Big(1-\frac{N}{1+s} 2^{1+s}\epsilon\Big).
\end{eqnarray*}
Thus we conclude that 
\begin{eqnarray*}
\epsilon^{-s}\cG_{s,\infty} (x,\epsilon e_1)>\frac{k_{N,s}}{2}4^s s^{-1}    \big(1-2^{1+s}\epsilon\big) |x|^{-N} x_1^s
 \end{eqnarray*}
 and for some $c>0$
 \begin{eqnarray*}
\epsilon^{-s}\cG_{s,\infty} (x,\epsilon e_1)<\frac{k_{N,s}}{2} 4^s s^{-1}  \big(1+c\epsilon^s |x|^{-s}  \big)     |x|^{-N}x_1^s,
 \end{eqnarray*}
 which imply that  for  $\epsilon>0$ small
 $$\Big| \epsilon^{-s}\cG_{s,\infty} (x,\epsilon e_1)-\cP_s(x)\Big| \leq   c\epsilon^s     |x|^{-N}x_1^s\leq c\epsilon^s   |x|^{s-N}, $$
 where $c>0$ is independent of $R$.
Thus,  (\ref{e 2.0-+}) holds uniformly in any compact set of $\R^N_+$.

 Now for any $R>1>\sigma>0$ large enough 
  \begin{eqnarray*}
\int_{B_{R}^+\setminus B_{\sigma} } \Big|\epsilon^{-s}\cG_{s,\infty} (x,\epsilon e_1)-\cP_s(x) \Big|dx &<&c\epsilon^s \int_{B_{R}^+\setminus B_{\sigma} }    |x|^{s-N}  dx
\\[2mm]&\leq &c   \epsilon^s  ( R^s-\sigma^s)
 \end{eqnarray*}
 and
 \begin{eqnarray*}
\int_{  B_{\sigma}^+ } \Big|\epsilon^{-s}\cG_{s,\infty} (x,\epsilon e_1)-\cP_s(x) \Big|dx &\leq &\int_{  B_{\sigma}^+ } \Big(\epsilon^{-s}\cG_{s,\infty} (x,\epsilon e_1)+\cP_s(x) \Big)dx
\\[2mm]&\leq &c \big( \epsilon^{-s}\sigma^{2s}  + \sigma^s\big).
 \end{eqnarray*}
Taking $\sigma=4\epsilon $ and $R=\epsilon^{-\frac23}$ with $\epsilon\in(0,\frac18)$, we can see that 
   \begin{eqnarray*}
\int_{B_{R}^+  } \Big|\epsilon^{-s}\cG_{s,\infty} (x,\epsilon e_1)-\cP_s(x) \Big|dx &\leq & c   \Big( \epsilon^s    R^s +\epsilon^{-s}\sigma^{2s}  + \sigma \Big)
\\[2mm]&\leq &c   \Big( \epsilon^{\frac s2}  +\epsilon^{ s }     \Big)
\\[2mm]&\to &0\quad{\rm as}\ \, \epsilon\to0^+
 \end{eqnarray*}
 and
 \begin{eqnarray*}
&&\int_{\R^N_+\setminus B_{R}^+  } \Big|\epsilon^{-s}\cG_{s,\infty} (x,\epsilon e_1)-\cP_s(x) \Big|\frac{1}{(1+|x|)^{N+2s}} dx 
\\[2mm] &\leq & \int_{\R^N_+\setminus B_{R}^+  } \Big(\epsilon^{-s}\cG_{s,\infty} (x,\epsilon e_1)+\cP_s(x) \Big)\frac{1}{(1+|x|)^{N+2s}} dx
\\[2mm]&\leq &c   \int_{\R^N_+\setminus B_{R}^+  } \Big(\epsilon^{-s}|x-\epsilon e_1|^{2s-N}+|x|^{s-N} \Big)\frac{1}{(1+|x|)^{N+2s}} dx 
\\[2mm]&\leq & c\big( \epsilon^{-s} R^{2s-N}+R^{s-N}\big)R^{-2s}
\\[2mm]&\leq & c\big( \epsilon^{\frac23 N-s}   +\epsilon^{-\frac23(s+N)}\big) 
\\[2mm]&\to &0\quad{\rm as}\ \, \epsilon\to0^+.
 \end{eqnarray*}
 For  $\varphi\in \bX_s(\R^N_+)$, we have that $(-\Delta)^s\varphi $ is bounded in $\R^N_+$ and 
 $$|(-\Delta)^s\varphi(x) |\leq \frac{c}{(1+|x|)^{N+2s}}$$ for some $c>0$,  then we deduce (\ref{e 2.1-00}).

 Therefore,  (\ref{e 2.0-+}) holds in $L^1_s(\R^N_+)$. 
 We complete the proof. \hfill$\Box$\medskip
 
 In the approximation of the fundamental solution $\cP_s$, we need the following regularity results:
 
\begin{lemma}\label{lm 2.1-0} 
Assume that $w\in C^{2s+\epsilon}(\bar B_1)$ with $\epsilon>0$
 satisfies
 $$(-\Delta)^s w=f\quad {\rm in}\quad B_1,$$
 where $f\in C^1(\bar B_1)$. Then for $\beta\in (0,s)$, there exist  $ c_1  , c_2>0$ such that
\begin{equation}\label{2.0}
\|w\|_{C^\beta(\bar B_{1/2})}\le {c_1}\Big(\|w\|_{L^\infty(B_1)}+\|f\|_{L^\infty(B_1)}+\| w\|_{L^1_s(\R^N)}\Big)
\end{equation}
and
\begin{equation}\label{2.0-0}
\|w\|_{C^{2s+\beta}(\bar B_{1/4})}\le {c_2}\Big(\|w\|_{L^\infty(B_1)}+\|f\|_{C^\beta(B_\frac12)}+\| w\|_{L^1_s(\R^N)}\Big).
\end{equation}
\end{lemma}
 The proof is postponed in the appendix.   \hfill$\Box$

\medskip

Now we are in a position to prove Theorem \ref{th 1.1}. \medskip
 
 \noindent {\bf Proof of Theorem \ref{th 1.1}. } From (\ref {2.0}), for any compact set $\cO$ in $\R^N_+$, there exists $c>0$ such that 
$$
\|\epsilon^{-s}\cG_{s,\infty} (\cdot,\epsilon e_1)\|_{C^{2s+\beta}(\cO)}\le {c_2}  \| \epsilon^{-s}\cG_{s,\infty} (\cdot,\epsilon e_1)\|_{L^1(\R^N,d\mu_s)}.
$$
Therefore,  up to subsequence, there holds
$$
\epsilon^{-s}\cG_{s,\infty} (\cdot,\epsilon e_1)\to \cP_s\quad{\rm in} \ \, C^{2s+\beta}(\cO)\quad{\rm as}\ \, \epsilon\to0^+,
$$
 which, together with (\ref{e 2.0-+}) in $L^1(\R^N_+,d\mu_s)$, implies that for any $x\in\R^N_+$
 $$(-\Delta)^s \cP_s(x)=\lim_{\epsilon\to0^+} \epsilon^{-s}(-\Delta)^s\cG_{s,\infty} (x,\epsilon e_1)=0. $$
 Therefore, we obtain that 
 $$  
 (-\Delta)^s \cP_s(x)=0 \quad    {\rm in }\ \, \R^N_+,  
\qquad\ \ \cP_s(x)=0 \ \ \,   {\rm in }\ \,  \R^N_{*}\setminus\{0\}.
 $$
 
 From (\ref{e 2.1-00}), we have that 
   \begin{eqnarray*}
    \int_{\R^N_+} \cP_s(x) (-\Delta)^s\varphi(x) dx&= & \lim_{\epsilon\to0^+}\int_{\R^N_+}\epsilon^{-s}\cG_{s,\infty} (x, \epsilon e_1)(-\Delta)^s\varphi(x) dx
  \\[1mm]&= &\lim_{\epsilon\to0^+}   \epsilon^{-s} \varphi(\epsilon e_1)=\lim_{\epsilon\to0^+}    \frac{\varphi(\epsilon e_1)-\varphi(0)}{\epsilon^{s}}
  \\[1mm]&= &\frac{\partial^s }{\partial x_1^s} \varphi(0).
  \end{eqnarray*}
 We complete the proof. \hfill$\Box$ 
 
    \setcounter{equation}{0}
    \section{Approximation by the Poisson kernel}

    In this section, we approximate the fundamental solution $\cP_s$ by the Poisson kernel
    and we consider the solution of 
\begin{equation}\label{eq 4.1t-1}
  \arraycolsep=1pt\left\{
\begin{array}{lll}
 (-\Delta)^s u  =0\quad   &{\rm in}
\ \,  \R^N_+,\\[2mm]\phantom{---\, }
u=\epsilon^s\delta_{-\epsilon e_1}\quad &  {\rm in}\ \,  \R^N_{*} 
\end{array}\right.
\end{equation}
for $\epsilon\in(0,1)$, where we recall $\R^N_{*}=(-\infty,0]\times\R^{N-1}$.

\begin{theorem}\label{th 4.1}
Problem  (\ref{eq 4.1t-1})   has  a unique nonnegative  solution 
$$u_\epsilon(x)=\arraycolsep=1pt\left\{
\begin{array}{lll}
\cK_{s} \epsilon^{-s} x_1^s |x+\epsilon e_1|^{-N}\quad    &{\rm in}
\ \,  \R^N_+,\\[2mm]\phantom{\, }
 \epsilon^s\delta_{-\epsilon e_1}\quad &  {\rm in}\ \,  \R^N_{*}.
\end{array}\right. 
$$  Moreover, there holds
\begin{equation}\label{e 4.1}
\int_{\R^N_+} u_\epsilon(x)(-\Delta)^s\xi(x)  dx =\epsilon^s
\int_{\R^N_+}\xi(x)\Gamma_\epsilon(x) dx,\quad \forall \,\xi\in \bX_s(\R^N_+),
\end{equation}
where
\begin{equation}\label{Gamma}
  \Gamma_\epsilon(x)=\frac{c_{N,s}}{|x+\epsilon e_1|^{N+2s}},\quad\forall\, x\in\R^N_+.
\end{equation}
\end{theorem}

Before the proof of Theorem \ref{th 4.1},  we consider the  solution of
\begin{equation}\label{eq 4.1t}
  \arraycolsep=1pt\left\{
\begin{array}{lll}
 (-\Delta)^s u  =0\quad   &{\rm in}
\quad  B_r(re_1),\\[2mm]\phantom{---\, }
u=t^s\delta_{-te_1}\quad &  {\rm in}\quad \R^N\setminus B_r(re_1),
\end{array}\right.
\end{equation}
where $t\in(0,1)$.   
\begin{proposition}\label{pr 4.1}
Problem  (\ref{eq 4.1t})   has  a unique nonnegative  solution 
\begin{equation}\label{e 4.2-e}
u_{t,r}(x)=\cK_{s}\left(\frac{2rx_1 -|x|^2}{2tr+t^{2}} \right)^s|x+t e_1|^{-N},\quad \forall x\in B_r(re_1)
\end{equation} 
and 
\begin{equation}\label{e 4.2}
\int_{B_r(re_1)} u_{t,r}(x)(-\Delta)^s \xi(x)  dx =t^{s}
\int_{B_r(re_1)}\xi(x)\Gamma_t(x) dx,\quad \forall \xi\in C^\infty_0(B_r(re_1)).
\end{equation} 
\end{proposition}
\noindent{\bf Proof. } {\it Step 1: approximation of the Dirac mass. }  Let $g_0:\R^N\to[0,1]$ be a radially symmetric decreasing $C^2$ function  with the support in $\overline{B_{\frac12}(0)}$
such that $\int_{\R^N} g_0(x)dx=1$.
For any $n\in\N$ and $t\in (0,1)$, we denote
$$g_n(x)=n^{N}g_0(n(x+te_1)),\qquad\forall x\in \R^N.$$
Then we certainly have  that
$$
g_n\rightharpoonup\delta_{-te_1}\ \ {\rm as}\ \ n\to+\infty
$$
in the distribution sense and for any $t>0$, there exists $m_t>0$ such that for any $n\ge m_t$,
$${\rm supp} (g_n)\subset \overline{B_\frac t2(-te_1)}.$$

For $t>0$,  problem
$$
  \arraycolsep=1pt\left\{
\begin{array}{lll}
 (-\Delta)^s u   =0\quad   &{\rm in}
\quad  B_r(re_1),\\[2mm]\phantom{---\, }
u=t^sg_n\quad &  {\rm in}\quad \R^N\setminus B_r(re_1) 
\end{array}\right.
$$
 admits a unique solution $w_n$.

Denote that
$$
\tilde g_n(x):=c_{N,s}   \int_{\R^N}\frac{g_n(y)}{|x-y|^{N+2s}}dy,\quad \forall x\in B_r(re_1).
$$
For $n\ge m_t$, we have that supp$(g_n)\subset B_{\frac t2}(-te_1)$ and then $\tilde g_n\in C^1(\overline{B_1(e_1)})$ and
  $$\tilde w_n=w_n-t^s g_n\quad {\rm in}\quad \R^N.$$
  By the definition of the fractional Laplacian, it implies that
\begin{eqnarray*}
  (-\Delta)^s \tilde w_n(x) &=&(-\Delta)^s w_n(x)-t^s (-\Delta)^s g_n(x)  \\
   &=& c_{N,s}t^s \int_{\R^N}\frac{g_n(z)}{|z-x|^{N+2s}} dz
   \\&=&t^s \tilde g_n(x).
\end{eqnarray*}

Then   $\tilde w_n$ is the unique solution of
$$
  \arraycolsep=1pt\left\{
\begin{array}{lll}
 (-\Delta)^s u   =t^s \tilde g_n\quad   &{\rm in}
\quad  B_r(re_1),\\[2mm]\phantom{--- \, }
u=0\quad &  {\rm in}\quad \R^N\setminus B_r(re_1)
\end{array}\right.
$$
and
\begin{equation}\label{e 4.2-0}
\int_{B_r(e_1)} \tilde w_n(x)(-\Delta)^s \xi(x)  dx =t^s
\int_{B_r(re_1)}\xi(x)\tilde g_n(x) dx,\quad \forall \xi\in C^\infty_0(B_r(re_1)).
\end{equation}

{\it Step 2:  we prove that $\tilde g_n$ converges to $\Gamma_t$
uniformly in $B_r(re_1)$ and in $C^{\theta}(B_r(re_1))$ for $\theta\in(0,1)$. }
 
   It is obvious that $\tilde g_n$ converges to $\Gamma_s$ every point in $ \overline{B_1(e_1)}$.
For $x,y\in B_r(re_1)$ and any $n\in\N$, we have that
\begin{eqnarray*}
&& |\tilde g_n(x)-\tilde g_n(y)| 
\\ &=&c_{N,s}|\int_{B_{\frac t2}(-te_1)}\Big[\frac{1}{|x-z|^{N+2s}}-\frac{1}{|y-z|^{N+2s}}\Big]g_n(z)dz|
 \\ &\le & c_{N,s}\int_{B_{\frac t2}(-te_1)}\frac{\big||x-z|^{N+2s}- |y-z|^{N+2s}\big|}{|x-z|^{N+2s}|y-z|^{N+2s}} g_n(z)dz
 \\ &\le & c_{N,s}(N+2s)|x-y|\int_{B_{\frac t2}(-te_1)}\frac{|x-z|^{N+2s-1}+ |y-z|^{N+2s-1}}{|x-z|^{N+2s}|y-z|^{N+2s}} g_n(z)dz
 \\ &\le & c_3|x-y|\int_{B_{\frac t2}(-te_1)}g_n(z)dz
 \\ &= & c_3|x-y|,
\end{eqnarray*}
where $c_3>0$ independent of $n$. So $\{\tilde g_n\}_n$ is uniformly bounded in $C^{0,1}(B_r(re_1))$. Combining the converging
$$\tilde g_n\to \Gamma_t\ {\rm every\ point\ in} \  \overline{B_r(re_1)}.$$
We conclude that $\tilde g_n$ converges to $\Gamma_s$
uniformly in $B_r(re_1)$ and in $C^{\theta}(B_r(re_1))$ for $\theta\in(0,1)$.
 
 {\it Step 3: passing  to the limit. } We  denote $\mathcal{O}_i$ the open sets with $i=1,2,3$ such that
 $$\mathcal{O}_1\subset \bar\mathcal{O}_1\subset \mathcal{O}_2\subset \bar\mathcal{O}_2\subset\mathcal{O}_3\subset \bar\mathcal{O}_3 \subset \R^N_+. $$
 By  Lemma \ref{lm 2.1-0} , for $\beta\in(0,s)$,  there exist  $c ,c' >0$ independent of $n$ such that
\begin{eqnarray*}
\norm{w_n}_{C^{\beta}(\mathcal{O}_2)} \le c [\norm{w_n}_{L^1(B_1(e_N))}+\norm{\tilde g_n}_{L^{\infty}(\mathcal{O}_3)}+\norm{w_n}_{L^{\infty}(\mathcal{O}_3)}] 
  \le c' 
\end{eqnarray*}
and
\begin{eqnarray*}
 \norm{w_n}_{C^{2s+\beta}(\mathcal{O}_1)} \le  c [\norm{w_n}_{L^1(B_1(e_N))}+\norm{\tilde g_n}_{C^{\beta}(\mathcal{O}_2)}
+\norm{w_n}_{C^{\beta}(\mathcal{O}_2)}] 
 \le c'. 
\end{eqnarray*}

Therefore,  by the Arzela-Ascoli Theorem, there exist $u_{t,r}\in C^{2s+\epsilon}_{\rm loc}$ in $B_r(re_1)$ for some $\epsilon\in(0,\beta)$ and a subsequence $\{w_{n_k}\}$ such that
$$
w_{n_k}\to u_{t,r}\quad  {\rm \ in}\quad C^{2s+\epsilon}\ {\rm locally\ in }\ \R^N_+, \quad{\rm as}\quad  n_k\to\infty.
$$
Passing the limit of (\ref{e 4.2-0})  as $n_k\to\infty$, we obtain (\ref{e 4.2}).  

The Poisson kernel of $B_r$ (see \cite{BGR})   has the formula
\begin{equation}\label{eq:Poisson-ker-BR}
\cP_{s,r}(x,y)=\arraycolsep=1pt\left\{
\begin{array}{lll}
\cK_{s}\left(\frac{r^2-|x|^2}{|y|^2-r^2} \right)^s|x-y|^{-N}\quad    &{\rm if}
\ \,  |y|>r,\,|x|<r,\\[2mm]\phantom{\, }
0\quad &  {\rm if\ not}.
\end{array}\right. 
\end{equation}
  The constant $\cK_{s}$ is chosen such
that
$$
\int_{\R^N}\cP_{s,r}(0,y)\,dy = \int_{\R^N \setminus B_r}\cP_{s,r}(0,y)\,dy =
1.
$$
Then 
$$u_{t,r}(x)=\cP_{s,r}\big(x-re_1,-(t+r)e_1\big)=\cK_{s}\left(\frac{2rx_1 -|x|^2}{2tr+t^{2}} \right)^s\big|x+t  e_1\big|^{-N}. $$
 We complete the proof.   \hfill$\Box$\medskip

\noindent{\bf Proof of Theorem \ref{th 4.1}. } 
 Since 
 $$B_r(re_1)\subset B_R(Re_1)\ \ {\rm for}\ \,  R\geq r\quad{\rm and}\ \  \R^N_+=\bigcup_{r>1} B_r(re_1).$$
Taking $t=\epsilon$ in (\ref{e 4.2-e}), we have that 
 \begin{equation}\label{e 4.2-e}
u_{\epsilon,r}(x)=\cK_{s}\left(\frac{2x_1 -\frac{|x|^2}{r}}{2\epsilon+\frac{\epsilon^{2}}{r}} \right)^s|x+\epsilon  e_1|^{-N},\quad \forall x\in B_r(re_1),
\end{equation} 
then $r\to u_{\epsilon,r}(x)$ is nondecreasing for any $\epsilon$ and $x$. Let
\begin{equation}\label{eq 1.1-po0}
 u_\epsilon(x)=\arraycolsep=1pt\left\{
\begin{array}{lll}
\cK_{s} \epsilon^{-s} x_1^s\, |x+\epsilon e_1|^{-N}\quad    &{\rm in}
\ \,  \R^N_+,\\[2mm]\phantom{\, }
 \epsilon^s\delta_{-\epsilon e_1}\quad &  {\rm in}\ \,  \R^N_{*}
\end{array}\right. 
\end{equation}
and then  we note that the sequence $\{u_{\epsilon,r}\}_{r>0}$ has an upper bound  $u_\epsilon$ in
$\R^N_+$.

 By the direct computation,  as $r\to+\infty$,  $u_{t,r}$   converges to $u_\epsilon\in C^{2s+\theta}$ in any compact set of $\overline{\R^N_+}$ and in $L^1_s(\R^N_+)$. Therefore, $u_\epsilon$ is  a solution of  (\ref{eq 4.1t}).

Then for any $\varphi\in C_c^\infty(\R^N_+)$, there exists $r>0$ such that 
 the support of $\varphi$ is a subset of $B_r(re_1)$ and passing to the limit of (\ref{e 4.2}),
 we obtain that 
$$
\int_{\R^N_+} u_\epsilon(x)(-\Delta)^s\xi(x)  dx =\epsilon^s
\int_{\R^N_+}\xi(x)\Gamma_\epsilon(x) dx,\quad \forall \,\xi\in C^\infty_c(\R^N_+).
$$
Finally, this identity holds for  $\xi\in \bX_s(\R^N_+)$ since $C^\infty_c(\R^N_+)$ is dense in $\bX_s(\R^N_+)$.\hfill$\Box$\medskip

From the above proof, we have the following corollary.

\begin{corollary} \label{cr 3.1}
Let  $\mu $ be a Radon measure with the support in $\R^N_-$ and
$$\cP_{s,\infty}(x,y)=\cK_{s}\left( \frac{x_1}{-y_1}\right)^s|x-y|^{-N}\qquad{\rm for}\ \,  x\in\R^N_+,\, y\in \R^N_-.$$

Then    
$$ \cP_{s,\infty}[\mu ](x):=\left\{
 \begin{array}{lll}
 \int_{\R^N_*}\cP_{s,\infty}(x,y)   d\mu(y)\quad &  {\rm for }\ \, x\in\R^N_+, \\[2.5mm]
 \mu \quad &  {\rm in }\ \  \R^N_-
 \end{array}
 \right.
 $$ 
is the unique solution of 
\begin{equation}\label{eq 1.1-po}
 \left\{
 \begin{array}{lll}  
 (-\Delta)^s u=0 \quad &  {\rm in }\ \, \R^N_+,  \\[2.5mm] 
\qquad\ \ u=\mu  & {\rm in }\ \,   \R^N_{*}
\end{array}
 \right.
\end{equation}
and the following distributional identity holds
\begin{equation}\label{e 4.00}
\int_{\R^N_+}\cP_{s,\infty}[\mu] (x)(-\Delta)^s\xi(x)  dx = 
\int_{\R^N_+}\xi(x)\Gamma_\mu (x) dx,\quad \forall \,\xi\in \bX_s(\R^N_+),
\end{equation}
where 
$$\Gamma_\mu (x)=c_{N,s} \int_{\R^N_*}\frac{d\mu(y)}{|x-y|^{N+2s}}.$$
\end{corollary}

 \noindent {\bf Proof of Theorem \ref{th 1.2}. }  
   From Lemma \ref{lm 2.1} 
\begin{equation}\label{e 2.0}
\epsilon^{-s}\cG_{s,\infty} (\cdot, \epsilon e_1)\to \cP_s\quad  \  {\rm as}\ \, \epsilon\to0^+
\end{equation}
in $L^1_s(\R^N_+)$ and uniformly in any compact set of $\R^N_+$.

We next show 
  $$\epsilon^{s}\cP_{s,\infty}(\cdot, -\epsilon e_1)\to \cP_s\quad  \  {\rm as}\ \, \epsilon\to0^+ $$
in $L^1_s(\R^N_+)$ and uniformly in any compact set of $\R^N_+$.

Recall that $u_\epsilon$, defined in (\ref{eq 1.1-po0}), is  the solution of  (\ref{eq 4.1t}) and
$$u_\epsilon(x)=\epsilon^{s}\cP_{s,\infty} (x, -\epsilon e_1)= \cK_{s} x_1^s|x+\epsilon e_1|^{-N}, $$
then $\epsilon^{s}\cP_{s,\infty} (\cdot, -\epsilon e_1)\in L^1_s(\R^N_+)$ and
$$\epsilon^{s}\cP_{s,\infty} (x, -\epsilon e_1)-\cP_s(x)=\frac{x_1^{\frac{N}{2}}+\epsilon^{\frac{N}{2}} }{|x|^N|x+\epsilon e_1|^N}  x_1^s. $$
Thus,  for any open set $O$,  $\overline O\subset \R^N_+$, 
\begin{equation}\label{e 2.1-up}
\epsilon^{s}\cP_{s,\infty} (\cdot, -\epsilon e_1)\to  \cP_s \quad{\rm as}\ \ \epsilon\to0^+\quad {\rm uniformaly\ in}\ O.  
\end{equation}

Moreover, let  $\sigma=\epsilon^{\theta}$ with $\theta\in (0,\frac{N}{3N-2s})$ and we see that 
   \begin{eqnarray*}
    &&\int_{\R^N_+} \Big| x_1^s|x+\epsilon e_1|^{-N}-x_1^s|x|^{-N}\Big| (1+|x|)^{-N-2s}dx \\[2mm]&\leq &4^N\epsilon^{\frac{N}{2}} \int_{\R^N_+}  \frac{x_1^{\frac{N}{2}}+\epsilon^{\frac{N}{2}} }{|x|^N|x+\epsilon e_1|^N}  x_1^s (1+|x|)^{-N-2s} dx  
    \\[2mm]&\leq&4^N\epsilon^{\frac{N}{2}}\Big( \int_{  B_\sigma^+(0)} \frac{x_1^{\frac{N}{2}}+\epsilon^{\frac{N}{2}}}{|x|^{N-\epsilon}|x+\epsilon e_1|^N}  dx+  
    \int_{\R^N_+\setminus B_\sigma(0)}  \frac{ (1+|x|)^{-N-2s} }{|x|^{N-s}|x+\epsilon e_1|^{\frac{N}2}}   dx \Big)
    \\[2mm]&\leq&4^N\epsilon^{\frac{N}{2}}\Big( c\sigma^{s-\frac{N}{2}} +\sigma^{s-\frac{3N}{2}}
    \int_{\R^N}  (1+|x|)^{-N-2s}dx\Big)
    \\[2mm]&\leq& c\epsilon^{\frac{N}{2}-(s-\frac{3N}{2})\theta}
    \\[2mm]&\to &0\quad{\rm as}\ \, \epsilon\to0^+.
\end{eqnarray*}

Moreover,   we recall that 
\begin{equation}\label{e 4.1--}
\int_{\R^N_+} u_\epsilon(x)(-\Delta)^s\xi(x)  dx =\epsilon^s
\int_{\R^N_+}\xi(x)\Gamma_\epsilon(x) dx,\quad \forall \,\xi\in \bX_s(\R^N_+).
\end{equation}
The left hand side of (\ref{e 4.1--}) verifies the convergence  
$$\int_{\R^N_+} u_\epsilon(x)(-\Delta)^s\xi(x)  dx\to \int_{\R^N_+} \cP_s (x)(-\Delta)^s\xi(x)  dx\quad{\rm as}\ \, \epsilon\to0^+$$
and the right hand has
$$\epsilon^s
\int_{\R^N_+}\xi(x)\Gamma_\epsilon(x) dx\to \frac{\partial^s}{\partial x_1^s} \xi(0)\quad{\rm as}\ \, \epsilon\to0^+,$$
then we deduce the identity (\ref{id 1.1}) from the Poisson kernel.
  \hfill$\Box$

     \setcounter{equation}{0}
    \section{The solution  blowing up on whole boundary }
  \subsection{Distributional indentity}
  
  \noindent{\bf Proof of Theorem \ref{th 1.3}. }
 Let 
 $$ \mu_t=\delta_{-t}(y_1) \big(t^s \omega_{\R^{N-1}}(y')\big) \quad{\rm with}\ \, t>0,$$
 where $\omega_{\R^{N-1}}(y')$ is the Hausdorff measure of $\R^{N-1}$ and $d\omega_{\R^{N-1}}(y')=dy'$ and $\delta_t(h)=h(t)$ for $h\in C_c(\R)$.  From the Poisson kernel expression, we obtain that 
    \begin{eqnarray*}
    \cQ_{s,t}(x_1,x')&=& \int_{\R^N_*}\cP_{s,\infty}(x,y) d\mu_t
     \\[2mm]&=&  \cK_{s} x_1^s \int_{\R^{N-1}}    \big((x_1+t)^2+|x'-y'|^2\big)^{-\frac{N}2} dy'
     \\[2mm]&=&\cC_1 x_1^{s} (x_1+t)^{-1}  
     \\[2mm]&\to &\cC_1x_1^{s-1}  \quad {\rm as}\ \, t\to 0^+
 \end{eqnarray*}
 and
the above convergence holds in $L^1_s(\R^N_+)$ and uniformly in any compact set of $\R^N_+$,
   where
  $$\cC_1=\cK_{s} \int_{\R^{N-1}} \big(1+|z'|^2\big)^{-\frac{N}2} dz'. $$
  Recall that 
  \begin{eqnarray*}
  \Gamma_{f_1}(x)&=&c_{N,s} \int_{\R^N}\frac{1}{|x-y|^{N+2s}}d\mu_t(y)
  \\[2mm]&=&c_{N,s} \int_{\R^{N-1}}\frac{t^s}{\big((x_1+t)^2+|x'-y'|^2\big)^{\frac{N+2s}2}}dy'
    = \cC_2 t^s (x_1+t)^{-2s}, 
  \end{eqnarray*}
  where 
  $$\cC_2=c_{N,s}\int_{\R^{N-1}} \big(1+|z'|^2\big)^{-\frac{N+2s}2} dz'.$$
Then (\ref{e 4.00}) with $\mu = \mu_t$, implies that   for any $\xi\in \bX_s(\R^N_+)$, 
     \begin{eqnarray*}
     \int_{\R^N_+}\cQ_{s,t}(x_1,x')  (-\Delta)^s\xi(x)  dx =\cC_2  
t^s\int_{\R^N_+}\xi(x)(x_1+t)^{-2s}  dx. 
   \end{eqnarray*}
   Note that 
   $$\cC_2  
t^s \int_{\R^N_+}\xi(x)(x_1+t)^{-2s}  dx\to \cC_2\int_{\R^{N-1}} \frac{\partial^s}{\partial {x_1^s}}\xi(0,x')dx'. $$
   Therefore, we we obtain that  
   \begin{eqnarray*}
     \int_{\R^N_+}x_1^{s-1} (-\Delta)^s\xi(x)  dx =\frac{\cC_2}{\cC_1}\int_{\R^{N-1}} \frac{\partial^s}{\partial {x_1^s}}\xi(0,x')dx'\quad{\rm for}\ \, \xi\in \bX_s(\R^N_+),
   \end{eqnarray*}
   where
 \begin{eqnarray*} \frac{\cC_2}{\cC_1}&=&   \frac{c_{N,s}\int_{\R^{N-1}} \big(1+|z'|^2\big)^{-\frac{N+2s}2} dz'}{\cK_s\int_{\R^{N-1}} \big(1+|z'|^2\big)^{-\frac{N}2} dz'}
  \\[2mm] &=&    \frac{c_{N,s}\omega_{_{N-2}} \cB(\frac{N+2s}{2}-\frac{N-1}{2}, \frac{N-1}{2}) }{\cK_s\omega_{_{N-2}} \cB(\frac{N}{2}-\frac{N-1}{2}, \frac{N-1}{2}) } 
 = 2 \sqrt{\pi} \, \frac{s}{\sin(\pi s)}\frac{\Gamma(s+\frac12)}{\Gamma(s)}, 
   \end{eqnarray*}
  here $\cB$ is the beta function and we use the fact that for $M\in\N$ and $\tau<-\frac{M}{2}$,
  \begin{eqnarray*}  
  \int_{\R^{M}} \big(1+|z|^2\big)^\tau dz&=&\omega_{_M} \int_0^{+\infty} (1+r^2)^\tau r^{M-1} dr    
  \\[2mm] &=&  \omega_{_M}   \int_0^1 t^{-\tau-1-\frac{M}{2}} t^{\frac{M}{2}-1} dr   
 = \omega_{_M} \cB(-\tau-\frac{M}{2}, \frac{M}{2})
   \\[2mm] &&\qquad \qquad\qquad\qquad\qquad\quad =  \omega_{_M} \frac{\Gamma(-\tau-\frac{M}{2})\Gamma(\frac{M}{2})}{\Gamma(-\tau)},
   \end{eqnarray*}
   where $\omega_{_M}$ is the volume of the unit sphere in $\R^M$.  \hfill$\Box$\medskip
  
   Let $$\cR_s(x)=(x_1)_+^s$$
    and it is known that $\cR_s$ is a solution of 
\begin{equation}\label{eqq ex}
(-\Delta)^s u=0 \quad    {\rm in }\ \,  \R^N_+,   
\qquad u=0  \ \   {\rm in }\ \,  \R^N_*.  
\end{equation}
Next we use our approximation to show this fact. 

  \begin{corollary}\label{th 1.4}
The function   $\cR_s$ is a solution of  (\ref{eqq ex})
 and it verifies that  
\begin{equation}\label{id 1.3} 
\int_{\R^N_+} \cR_s(x) (-\Delta)^s\varphi(x) dx=0\quad{\rm for\ any}\ \,  \varphi\in \bX_s(\R^N_+).
\end{equation}
  \end{corollary}

   \noindent {\bf Proof. } We take 
   $$ \nu_t=\delta_{-t}(y_1) \,  \big( t^{1+s}\omega_{\R^{N-1}}(y')\Big)  \quad{\rm with}\ \, t>0.$$  
   From the Poisson kernel expression, we obtain that 
    \begin{eqnarray*}
  \int_{\R^N_*}\cP_{s,\infty}(x,y)  d\nu_t
     &=&  \cK_{s}   x_1^{s }t \int_{\R^{N-1}}    \big((x_1+t)^2+|x'-y'|^2\big)^{-\frac{N}2} dx'
     \\[2mm]&=&\cC_1 x_1^{s } t (x_1+t)^{-1}  
     \\[2mm]&\to &\cC_1x_1^{s}  \quad {\rm as}\ \, t\to 0^+
 \end{eqnarray*}
 and 
the above convergence holds in $L^1_s(\R^N_+)$ and uniformly in any compact set of $\R^N_+$,
   where
  $$\cC_1=\cK_{s} \int_{\R^{N-1}} \big(1+|z'|^2\big)^{-\frac{N}2} dz'. $$
  Recall that for $x\in \R^N_+$
  \begin{eqnarray*}
  \Gamma_{\nu_t}(x) =  c_{N,s} \int_{\R^{N-1}}\frac{ t^{1+s}}{\big((x_1+t)^2+|x'-y'|^2\big)^{\frac{N+2s}2}}dy'
    = \cC_1 t^{ 1+s}  (x_1+t)^{-2s}.  
  \end{eqnarray*}
From (\ref{e 4.00}) with $\mu = \nu_t$,  we obtain that    for any $\xi\in \bX_s(\R^N_+)$, 
     \begin{eqnarray*}
     \int_{\R^N_+}\cQ_{s,t}(x_1,x')  (-\Delta)^s\xi(x)\,  dx =\cC_2  t^{1+s}
\int_{\R^N_+}\xi(x)  (x_1+t)^{-2s}  dx. 
   \end{eqnarray*}
   Note that 
   $$\cC_2  t^{1+s }
\int_{\R^N_+}\xi(x) (x_1+t)^{-2s}  dx\to 0\quad{\rm as}\ \, t\to0^+.  $$
   Therefore, we we obtain that for any $\xi\in \bX_s(\R^N_+)$
   \begin{eqnarray*}
     \int_{\R^N_+}x_1^{s} (-\Delta)^s\xi(x)  dx =0.
   \end{eqnarray*}
 We complete the proof.  \hfill$\Box$

 \subsection{ More $s$-harmonic functions in $\R^N_+$}
 
 It is shown in \cite[Theorem 5.1]{CLO} that the nonnegative solution of (\ref{eqq ex})
 only has the form 
 $$u=c\cR_s\quad {\rm with}\ c\geq0.$$
 In this subsection, we will show more $s$-harmonic functions in $\R^N_+$ without the restriction of positivity.  To this end,  we consider the functions with the separable  variables with the  form $u(x_1,x')=x_1^s h(x')$ and we can take an equivalent definition of the fractional Laplacian in the  principle value sense
  \begin{equation}\label{def 3.1} 
  (-\Delta)^s   u(x_1,x')=c_{N,s} \lim_{\epsilon\to0^+} \int_{(-\frac{1}{\epsilon},-\epsilon)\cup 
  (\epsilon,\frac{1}{\epsilon})} \int_{B^*_\frac1\epsilon  \setminus B^*_\epsilon  }\frac{ u(x_1,x')-
u(x_1+z_1,x'+z')}{|z|^{N+2s}}  dz'dz_1,
\end{equation}
where $B^*_r$ is the ball centered at the origin with radius $r$ in $\R^{N-1}$. \medskip

 \noindent{\bf Proof of Theorem \ref{th 1.2-p}. } 
 Let  $h_0 \in\cH^s(\R^{N-1})$ with $N\geq3$,
 $$u_0(x_1,x')= (x_1)_+^s h_0(x') $$
and we need prove that
 $$ (-\Delta)^s   u_0 =0\quad {\rm in}\ \ \R^N_+.   $$
 We use the definition (\ref{def 3.1}) to obtain that
  \begin{eqnarray*}
&&\frac1{c_{N,s}} (-\Delta)^s   u_0(x_1,x') \\[2mm] & =& \lim_{\epsilon\to0^+} \int_{(-\frac{1}{\epsilon},-\epsilon)\cup 
  (\epsilon,\frac{1}{\epsilon})} \Big((x_1)_+^s-(x_1+z_1)_+^s\Big)  \int_{B^*_\frac1\epsilon  \setminus B^*_\epsilon  }\frac{h_0(x')}{(|z_1|^2+|z'|^2)^{\frac{N+2s}2}}  dz'dz_1
     \\[2mm]&&+\lim_{\epsilon\to0^+} \int_{(-\frac{1}{\epsilon},-\epsilon)\cup 
  (\epsilon,\frac{1}{\epsilon})}  (x_1-z_1)_+^s   \int_{B^*_\frac1\epsilon  \setminus B^*_\epsilon  }\frac{h_0(x')-h_0(x'+z')}{(|z_1|^2+|z'|^2)^{\frac{N+2s}2}}  dz'dz_1 
   \\[2mm]&=&h_0(x') \lim_{\epsilon\to0^+} \int_{(-\frac{1}{\epsilon},-\epsilon)\cup 
  (\epsilon,\frac{1}{\epsilon})} \frac{(x_1)_+^s-(x_1+z_1)_+^s}{z_1^{1+2s}}dz_1  \int_{B^*_\frac1\epsilon  \setminus B^*_\epsilon  }\frac{1}{(1+|z'|^2)^{\frac{N+2s}2}}  dz' 
     \\[2mm]&&+\lim_{\epsilon\to0^+} \int_{(-\frac{1}{\epsilon},-\epsilon)\cup 
  (\epsilon,\frac{1}{\epsilon})}  (x_1-z_1)_+^s   \int_{B^*_\frac1\epsilon  \setminus B^*_\epsilon  }\Big(h_0(x')-h_0(x'+z')\Big) K_{z_1}(|z'|)  dz'dz_1.
   \end{eqnarray*}  
Note that 
$$ \int_{(-\frac{1}{\epsilon},-\epsilon)\cup 
  (\epsilon,\frac{1}{\epsilon})} \frac{(x_1)_+^s-(x_1+z_1)_+^s}{z_1^{1+2s}}dz_1 \to0\quad{\rm as}\ \, \epsilon\to0^+  $$
  and
 $$\int_{B^*_\frac1\epsilon  \setminus B^*_\epsilon  }\frac{1}{(1+|z'|^2)^{\frac{N+2s}2}}  dz'\to \int_{\R^{N-1}}\frac{1}{(1+|z'|^2)^{\frac{N+2s}2}}  dz'\quad{\rm as}\ \, \epsilon\to0^+.  $$
For any $z_1$, $K_{z_1}$ verifies (\ref{con 2.2-p}) in $\R^{N-1}$, see Remark \ref{rem 33-1} part $(i)$,  and from the proof of Proposition \ref{th 2.1-p}, we have that 
$$\int_{B^*_\frac1\epsilon  \setminus B^*_\epsilon  }\Big(h_0(x')-h_0(x'+z')\Big) K_{z_1}(|z'|)  dz'=0\quad{\rm for}\ \, \epsilon\in(0,1).$$
 Then passing to the limit, we obtain that 
 $$ (-\Delta)^s   u_0 =0\quad{\rm in}\ \ \R^N_+.$$

Set   
 $$v_0(x_1,x')= (x_1)_+^{s-1} h_0(x'),$$
via a similar proof, we can obtain that 
 $$ (-\Delta)^s   v_0 =0\quad {\rm in}\ \ \R^N_+.   $$
 
 We obtain that 
 $\big\{(x_1)_+^s v(x') :\ v \in \cH^s(\R^{N-1})\big\}\cup \big\{(x_1)_+^{s-1} v (x') :\ v \in \cH^s(\R^{N-1})\big\}\subset \cH^s(\R^N_+)$. We complete the proof.\hfill$\Box$\medskip

  \appendix
 \section{Appendix: Proof of Theorem \ref{lm 2.1}}

 The regularities of fractional Poisson problems have been studied in   \cite{CS-1,S,RS}.\medskip

     \begin{lemma}\label{app lm 1} 
 Let $f\in L^\infty (B_1)$ and $u:\R^N\to \R$ be a bounded function verifying
 $$(-\Delta)^s u=f\quad{\rm in}\ \, B_1.$$ 
 Then for $\beta\in(0,s)$
 $$\|u\|_{C^{\beta}(B_\frac12)}^*\leq C\big(\|u\|_{L^\infty(\R^N)}+\|f\|_{L^\infty(B_1)}\big). $$
 \end{lemma}
{\bf Proof.} It follows by  \cite[Theorem 12.1]{CS-1} replacing the   maximal and a minimal operator by the fractional Laplacian.  \hfill$\Box$  \medskip

 For $k\in\N$ and $\beta\in(0,1)$,  denote  
 $$\|u\|_{C^{k+\beta}(\Omega)}^*=\sum_{|\gamma|=k} \sup_{x,y\in\Omega} d^{k+\beta}_{x,y} \frac{D^\gamma u(x)-D^\gamma u(y)}{|x-y|^\beta}$$
 and for $b>0$
  $$\|u\|_{C^{\beta}(\Omega)}^{(b)}= \sup_{x\in\Omega} d^{b}_{x}|u(x)| + 
 \sup_{x,y\in\Omega} d^{b+\beta}_{x,y} \frac{D^\gamma u(x)-D^\gamma u(y)}{|x-y|^\beta}, $$
 where $\Omega$ is a bounded $C^{k+1}$ domain and $d^{k+\beta}_{x,y}=\min\{\rho(x),\rho(y)\}$.

 \begin{lemma}\label{app lm 2}\cite[Theorem 12.2.1]{CLO} (also see \cite{RS})
 Let $f\in C^\beta(\Omega)$ and $u:\R^N\to \R$ be a function verifying
 $$(-\Delta)^s u=f\quad{\rm in}\ \, \Omega.$$ 
 Then 
 $$\|u\|_{C^{2s+\beta}(\Omega)}^*\leq C\big(\|u\|_{L^\infty(\R^N)}+\|f\|^{(2s)}_{C^\beta(\Omega)}\big). $$

 \end{lemma}
 
 Now we are in a position to prove Theorem \ref{lm 2.1}. 
 \medskip
    
\noindent{\bf Proof of Theorem \ref{lm 2.1}. }     We denote $v=w\eta$, where $\eta:\R^N\to[0,1]$ is a $C^\infty$  function such that
$$\eta=1\quad {\rm in}\quad  B_{\frac12}\quad {\rm and}\quad   \eta=0\quad {\rm in}\quad  B_{\frac34}^c. $$
Then $v\in C^{2s+\epsilon}(\R^N)$ and for any $x\in B_{\frac12}$, $\epsilon\in(0,\frac18)$,
\begin{eqnarray*}
  (-\Delta)_\epsilon^s v(x)  = (-\Delta)_\epsilon^s w(x)+c_{N,s}\int_{\R^N \setminus B_\epsilon}\frac{\big(1-\eta(x+y)\big)w(x+y)}{|y|^{N+2s}}dy.
 \end{eqnarray*}
Together with the fact of $\eta(x+y)=1$ for $y\in B_\epsilon$, we derive that
 $$ \int_{\R^N \setminus B_\epsilon}\frac{\big(1-\eta(x+y)\big)w(x+y)}{|y|^{N+2s}}dy=\int_{\R^N  }\frac{\big(1-\eta(x+y)\big)w(x+y)}{|y|^{N+2s}}dy=: h_{1}(x),$$
thus, $$(-\Delta)^s v=h+c_{N,s} h_1\quad{\rm in}\quad B_{\frac12}.$$

$(i)$ For $x\in B_{\frac14}$ and $z\in\R^N\setminus B_{\frac12}$, there holds
\begin{eqnarray*}
  |z-x|\ge |z|-|x| \ge |z|-\frac14
   \ge \frac1{12}(1+|z|),
\end{eqnarray*}
which implies that
\begin{eqnarray*}
  |h_1(x)| = \Big|\int_{\R^N}\frac{(1-\eta(z))w(z)}{|z-x|^{N+2s}}dz\Big|
   &\le& \int_{\R^N\setminus B_{\frac12}}\frac{|w(z)|}{|z-x|^{N+2s }}dz \\
   &\le & 12^{N+2\alpha} \int_{\R^N}\frac{|w(z)|}{(1+|z|)^{N+2s}}dz\\
   &=&12^{N+2\alpha}\|w\|_{L^1_s(\R^N)}.
\end{eqnarray*}
From Lemma \ref{app lm 1},    there exists ${c_8}>0$ such that
\begin{eqnarray*}
\|v\|_{C^\beta(\bar B_{1/4})}&\le& {c_8}(\|v\|_{L^\infty(\R^N)}+\|h+h_1\|_{L^\infty(B_{1/2})}) \\[2mm]
   &\le& {c_8}(\|w\|_{L^\infty(B_1)}+\|h\|_{L^\infty(B_1)}+\|h_1\|_{L^\infty(B_{1/2})})\\[2mm]
   &\le& {c_9}(\|w\|_{L^\infty(B_1)}+\|h\|_{L^\infty(B_1)}+\| w\|_{L^1_s(\R^N)}),
\end{eqnarray*}
where  $ {c_9}=12^{N+2\alpha}{c_8}$.
Combining with $w=v$ in $B_{\frac12}$,  we obtain (\ref{2.0}).\smallskip

$(ii)$ For $x,\, y\in B_{\frac1{8}}$   and $|z|>\frac14$,  there exists $c>0$ such that 
  $$\frac1{|z-x|^{N+2s}}\leq \frac c{|z|^{N+2s}},\qquad\frac1{|z-y|^{N+2s}}\leq \frac c{|z|^{N+2s}} $$
  and
  $$  |z-x|^{N+2s}-|z-y|^{N+2s}\leq c|z|^{N-2s-1} |x-y|,$$
thus,
\begin{eqnarray*}
  |h_1(x)-h_1(y)| &=& \Big|  \int_{\R^N  }\frac{ \eta(x+z)w(x+z)-\eta(y+z)w(y+z)}{|z|^{N+2s}}dz  \Big|
  \\[2mm]
   &=&\Big|  \int_{\R^N\setminus B_\frac14  }\eta(z)w(z)\Big(\frac{1}{|z-x|^{N+2s}}-\frac{1}{|z-y|^{N+2s}}\Big) dz  \Big|
    \\[2mm]
   &\le & \int_{\R^N \setminus B_\frac14 } \big|  w(z) \big|  \frac{\big| |z-x|^{N+2s}-|z-y|^{N+2s}\big| }{|z-x|^{N+2s}|z-x|^{N+2s}}  dz      
     \\[2mm] &\leq & c\int_{\R^N \setminus B_\frac14 } \big|  w(z) \big|  \frac{ |x-y| }{1+|z|^{N+2s+1}}  dz
     \\[2mm] &\leq & c\|w\|_{L^1_s(\R^N)}|x-y|.
\end{eqnarray*}
Now we apply Lemma \ref{app lm 2} with $\Omega=B_{\frac1{8}}$ to obtain that 
\begin{eqnarray*} 
\|v\|_{C^{2s+\beta}(B_{\frac1{16}})} \leq c\|v\|_{C^{2s+\beta}(B_{\frac1{8}})}^*&\leq& C\big(\|v\|_{L^\infty(\R^N)}+\|h+h_1\|^{(2s)}_{C^\beta(B_\frac1{8})}\big)
\\[2mm]&\leq &c\big(\|w\|_{L^\infty(B_1)}+\|h\|_{C^\beta(B_\frac1{16})}+\| w\|_{L^1_s(\R^N)}\big).  
\end{eqnarray*}
which implies that 
\begin{eqnarray*} 
\|w\|_{C^{2s+\beta}(B_{\frac1{16}})} \leq  c\big(\|w\|_{L^\infty(B_1)}+\|h\|_{C^\beta(B_\frac1{16})}+\| w\|_{L^1_s(\R^N)}\big) 
\end{eqnarray*}
by the fact that $w=v$ in $B_{\frac12}$. Then  (\ref{2.0-0}) is proved. \hfill$\Box$

 \bigskip\bigskip
 
\noindent{\footnotesize  {\bf Acknowledgement: }  This work is supported by the Natural Science Foundation of China, No. 12071189, 12001252, by Jiangxi Province Science Funds, No. 20212ACB211005, 20202ACBL201001, by the Science and Technology Research Project of Jiangxi Provincial Department of Education, No. GJJ200307, GJJ200325.  
  }

\end{document}